\documentclass[twoside]{article}
\usepackage[a4paper]{geometry}
\usepackage[latin1]{inputenc} 
\usepackage[T1]{fontenc} 
\usepackage{RR}
\usepackage{hyperref}

\usepackage {times,amsmath,amssymb,latexsym,graphicx,epic,eepic,verbatim,stmaryrd}
\usepackage{multicol}
\usepackage{latexsym,color,verbatim}
\usepackage{amsmath}
\usepackage{amssymb}
\usepackage{amsfonts}
\usepackage{times}
\usepackage{graphicx,lscape}
\usepackage{epsfig}
\usepackage{array}
\usepackage{flafter}
\usepackage[all]{xy}

\newcommand{\expect}{\mathbb{E}}

\newcommand{\Preal}{\mathbb{R}^+}

\newcommand{\as}{\mbox{a.s.}}

\newcommand{\xii}{\mbox{\boldmath$\xi$}}

\newcommand\independent{\protect\mathpalette{\protect\independenT}{\perp}}
\def\independenT#1#2{\mathrel{\rlap{$#1#2$}\mkern2mu{#1#2}}}

\newtheorem{defn}{Definition}
\newtheorem{rem}{Remark}
\newtheorem{thm}{Theorem}
\newtheorem{prop}{Proposition}
\newtheorem{pf*}{Proof}
\newcommand{\qed}{$\square$}

\RRNo{8368}
\RRdate{September 2013}

\RRauthor{
Nicolas Tabareau\thanks{Inria, 
  Nantes, France}%
  \and
Jean-Jacques Slotine\thanks{Non Linear System Laboratory, Massachusetts Institute of
  Technology,
  Cambridge, Massachusetts, 02139, USA}}
\authorhead{Tabareau \& Slotine}
\RRtitle{Analyse de la contraction pour les systèmes dynamique
  alétoires non-linéaires}
\RRetitle{Contraction analysis of nonlinear random dynamical systems}
\titlehead{Contraction analysis of nonlinear random dynamical systems}
\RRresume{
  La théorie de la contraction (\cite{Lohmiller98} ) sert à l'étude de la stabilité des
  systèmes dynamiques non-linéaires.
  Dans ce rapport, nous étendons la théorie de la contraction
 aux cas des équations différentielles aléatoires. 
 Nous proposons deux nouvelles définitions de la contraction dans un
 cadre aléatoire (contraction presque sûre, et contraction aux
 moindres carrés) qui permettent de controler l'évolution d'un système
 stochastique de deux manières. La première garantie la convergence
 exponentielle pour presque toutes les réalisations du sysème, tandis
 que la deuxième garantie la convergence dans $L_2$ à une trajectoire
 unique. 
 Nous illustrons enfin la relative simplicité de cette extension en analysant
 des propriétés déterministes bien connu dans un cadre bruité. Plus
 spécifiquement, nous alnalysong la descente de gradient stochastique,
 l'impact du bruit sur la synchronisation d'oscillateurs et
 l'extension des propriétés de combinaison de systèmes contractant aux
 cas stochastique. 
}
\RRabstract{
In order to bring contraction analysis into the very fruitful and
topical fields of stochastic and Bayesian systems, we
extend here the theory describes in \cite{Lohmiller98} to random
differential equations.
We propose new definitions of contraction (almost sure contraction
and contraction in mean square) which allow to master the evolution of a
stochastic system in two manners. The first one guarantees eventual
exponential convergence of the system for almost all draws, whereas
the other guarantees the exponential convergence in $L_2$ of the
system to a unique trajectory.
We then illustrate the relative simplicity of this extension by analyzing usual
deterministic properties in the presence of noise. Specifically,
we analyze stochastic gradient descent, impact of
noise on oscillators synchronization and extensions of
combination properties of contracting systems to the stochastic
case. This is a first step towards combining the interesting and
simplifying properties of contracting systems with the probabilistic approach.
}
\RRmotcle{équations différentielles aléatoires, théorie de la contraction}
\RRkeyword{random differential equations, contraction theory}
\RRprojets{Ascola}
\RCRennes 

\begin{document}

\makeRR 

\tableofcontents

\section{Introduction}

The concept of random dynamical systems appeared to
be a useful tool for modelers of many different fields, especially in computational
neuroscience. Indeed, introducing noise in a system can 
be either a way to describe reality more sharply or to generate
new behaviors absent from the deterministic world. 
Naively, the definition of such systems generalized the usual
deterministic one as follows,
\begin{equation} \label{eq:random system}
{\bf x}_{i+1} = {\bf f}_i({\bf x}_i,i,\xii_i(\omega))
\qquad \mbox{or} \qquad
\dot{\bf x}_t = {\bf f}({\bf x}_t,\xii_t(\omega),t)
\end{equation} 
where $\xii_i$ (resp. $\xii_t$) defines a discrete
(resp. continuous) stochastic process.

Unfortunately, it comes out that whereas the extension of the
definition in the discrete time case is rather straightforward, things
get much harder for continuous time. Indeed, equation \ref{eq:random
  system} does not make sense for every continuous stochastic
process. There are basically to different ways of fixing this problem :
\begin{enumerate}
  \item
    One can restrict the definition to sufficiently smooth stochastic
    processes, namely those whose trajectories are \emph{cadlag} (see
    subsection \ref{sct:continuous}). This gives rise to \emph{random
      differential equations}.
  \item
    Or, as important stochastic processes like \emph{White noise} do
    not satisfy this property, one can switch to stochastic
    differential equations and It\^o calculus.
\end{enumerate}
When first looking at it, It\^o calculus seems more accurate as it just
generalized the usual differential calculus. But the problem of this
formalism lies in the absence of ``realisability'' of such defined
systems. By realisability, we mean the ability to simulate the
behavior of a system by mean of a computer. Hence, all result we can
obtain are purely theoretic, and it is hard to make up an intuition on
how such a system evolves. From a modeling point of view, this lake
of manageability make this formalism unusable in most fields as we are
hardly able to find explicit solution of a huge dynamical system as
encountered in biology, and thus even in the deterministic case. It
follows that we need our system to be realizable, and so
we must restrict our attention to random differential equations.

Once this restriction performed, we can really solve our system with
the traditional tools of differential analysis. Indeed, it now makes
sense to fix a draw in the probability space and solve the equation as
the noise being an external input. We can thus describes the
trajectory of the system given a sample path of the stochastic
process. It results that we are also able to use deterministic
contraction theory to study stability when the noise never put the
system out of contraction bounds. But this implies that every
trajectory is contracting, and so noise does not really matter. 

In this paper, we will investigate new definitions of stochastic
stability together with sufficient condition to guarantee that a
random differential system as a nice behavior even if the noise can
induce partial divergence in the trajectories. 

The notations of this paper follow those of the keystone
\cite{arnold98}.

\section{Part 1 : State's dependency}

\subsection{Nonlinear random system : the discrete time case}
\label{sct:discrete}

\subsubsection{Almost sure contraction : asymptotic exponential
  convergence $\as$} \label{sct:almost}

As a first step, we define the stochastic contraction in the field of
discrete system. In that case, there is no problem regarding the
equation generating the dynamic system as it is just the iterated
application of possibly different functions. We are dealing with
stochastic processes of the form :
$$
\xii : \Omega \times \mathbb{N} \rightarrow \mathbb{R}^n
$$
ie. we assume $\xii$ to be equivalent to a sequence of random variable
$\xii_i : \Omega \rightarrow \mathbb{R}^n$.

The definition of contraction in discrete-time case is a direct
extension of the deterministic case and makes use of the notion of
\emph{discrete stochastic differential systems} (\cite{jaswinski70})
of the form 
$${\bf x}_{i+1} = {\bf f}_i({\bf x}_i,i,\xii_i(\omega))$$
where the ${\bf f}_i$ are continuously differentiable (a condition
that will be assumed in the rest of this paper).

This definition is the natural extension of the definition of
\cite{Lohmiller98}, which assessed that the the difference between two
trajectories tends to zero exponentially fast. In the stochastic case,
we look for similar conditions satisfy almost surely. 

First, we have to reformulate the traditional property satisfied by
the metric allowing a space deformation. In the deterministic case,
the metric $M_i({\bf x}_i,i)$ has to be uniformly positive definite with
respect to ${\bf x}_i$ and $i$. This property of uniformity for a
metric $M_i({\bf x},i,\xii_i)$ depending also on noise can be written 
$$
\exists \lambda \quad \forall {\bf x},i \quad M({\bf x},i,{\bf
  \xii}_i) \geq \lambda {\bf I} \quad \as
$$

But as we want the noise to introduce local bad behaviors, we need to relax the
property in a sense that contraction can only been guaranteed
asymptotically. This introduces a slightly difficulty as the naive
formulation 
$$
\exists 0 \leq \alpha < 1 \qquad \lim_{n \rightarrow \infty}(\| \delta
{\bf z} \| -\| \delta {\bf z}_o \| \alpha^n) = 0
$$
just says that $\| \delta {\bf z} \|$ tends to zero, which is not what we
want. We thus have to switch to the following refined formulation.

\begin{defn}
  The random system ${\bf x}_{i+1} = {\bf f}_i({\bf
  x}_i,i,\xii_i(\omega))$, is said to be \emph{almost surely contracting}
  if there exists a uniformly positive definite metric $M_i({\bf
    x},i,\xii_i)$ and $\eta < 0$ such that: 
  $$  
  \mathbb{P}\{\omega \in \Omega, \ \limsup_{n \rightarrow \infty}
  \frac{1}{n} \log(\| \delta {\bf z} \|) \leq \eta \} = 1
  $$
  ie. the difference between two trajectories tends almost surely 
  to zero exponentially.
\end{defn}

\begin{rem}
\begin{enumerate}
\item
  The notion of contracting region cannot be extended to stochastic case
  as it is hardly possible to guarantee that a stochastic trajectory
  stay in a bowl without requiring strong bound on the noise. And in
  that case, the noise can be treated as a bounded perturbation by
  analyzing the worst case.
\item
  The definition makes use of logarithm whereas $\| \delta {\bf z} \|$ can
  be equal to $0$. Nevertheless, the reader should not be deterred and
  each time such a case appears, the equation is also satisfy if we allow
  infinite value and basic analytic extension (for example, $\expect(
  \log \alpha) =- \infty$ as soon as $\alpha(\omega) = 0$ for some
  $\omega$).
\end{enumerate}
\end{rem}

We can now state the first theorem of this paper. Remark that in the
definition of the system below, the dependence on the stochastic
perturbation is almost linear, as we can not master a malicious
non-linear system which strongly use the ``bad draw'' of the
stochastic process to diverge. In a sense, the system must satisfy a
notion of uniformity with respect to the stochastic process.

\begin{thm}\label{thm:ascDiscret}

  Given the random system ${\bf x}_{i+1} = {\bf f}_i({\bf
  x}_i,i,\xii_i(\omega))$,
  note $\sigma_{f}({\bf x}_i,i,\xii_i)$ the largest singular
  value of the generalized Jacobian of ${\bf f}_i$ at time i 
  according to some metric $M_i({\bf x},i,\xii_i)$. \\

  A sufficient condition for the system to be \emph{almost surely
    contracting} is that
  \begin{itemize}
  \item
    the random process $\sigma_{f_i}({\bf x},i,\xii_i)$ can
    be bounded independently from ${\bf x}$, ie there exists a
    stochastic process $\eta_i$ such that
    $$
    \forall {\bf x} \quad \sigma_{f_i}({\bf x},i,\xii_i) \leq
    \eta_i \qquad \as
    $$
  \item
    the stochastic process $\log(\eta_i)$ follows the strong
    law of large number (eg. i.i.d.)
    $$ \frac{1}{n} \sum_{i=1}^n \log(\eta_i)
    \rightarrow \frac{1}{n} \sum_{i=1}^n
    \expect(\log(\eta_i)) \qquad \as
    $$
  \item
    the expectation of the random variables $\log(\eta_i)$ can
    be uniformly bounded by some $\eta < 0$
    $$
    \forall i \quad \expect(\log(\eta_i)) \leq \eta
    $$
  \end{itemize}
  
\end{thm}

\begin{pf*}{Proof.}
  
  We make a strong use of the basic idea of the original proof.

  Note $F_i$ the discrete generalized Jacobian of $f_i$ :
  $ F_i({\bf x},i,\xii_i) \ = \ \Theta_{i+1}
  \frac{\partial}{\partial {\bf x}} {\bf f}_i({\bf x},i,{\bf
    \xii}_i) \Theta^{-1}_i$
  we have:
  $$
  \begin{array}{lrcl}
    & {\bf \delta {\bf z}}_{i+1}^T {\bf \delta {\bf z}}_{i+1} & = &
    {\bf \delta {\bf z}}_i^T (
    F_i^T F_i) \ {\bf \delta {\bf z}}_i \\
    \\
    \Rightarrow \quad & {\bf \delta {\bf z}}_{i+1}^T {\bf \delta {\bf z}}_{i+1} &
    \le & \sigma^2_{f_i}({\bf x}_i,i,\xii_i)
    \ {\bf \delta {\bf z}}_i^T {\bf \delta {\bf z}}_i
    \qquad \as
  \end{array}
  $$
  
  and hence,
  
  $$
  \| \delta {\bf z}_n \|^2 \le \| \delta {\bf z}_o \|^2 \prod_{i=0}^n
  \sigma^2_{f_i}({\bf x}_i,i,\xii_i)
  \qquad \as
  $$

  So by monotony of logarithm and the two required properties, we can
  deduce that for almost every $\omega$
  $$
  \begin{array}{rcl}
    \limsup_{n \rightarrow \infty} \frac{1}{n} \log(\| \delta {\bf
      z}_n \|) & \leq & \limsup_{n \rightarrow \infty}
    \frac{1}{n} \log(\prod_{i=0}^n \sigma_{f_i}({\bf x}_i,i,{\bf
      \xii}_i)) \\ 
    & \leq & \limsup_{n \rightarrow \infty}
    \frac{1}{n} \sum_{i=0}^n \log(\eta_i) \\ 
    & = & \limsup_{n \rightarrow \infty} \frac{1}{n} \sum_{i=0}^n 
    \expect(\log(\eta_i)) \\
    & \leq & \eta \\
  \end{array}
  $$
  That is
  $$
  \mathbb{P}\{\omega \in \Omega, \ \limsup_{n \rightarrow \infty}
  \frac{1}{n} \log(\| \delta {\bf z} \|) \leq \eta \} = 1
  $$
\qed
\end{pf*}

\subsubsection{Contraction in mean square: asymptotic exponential
  convergence in mean square}

We have seen in the subsection above sufficient conditions to guarantee
almost sure asymptotic exponential convergence. But we could also be
interested in looking for conditions to guarantee exponential
convergence in mean square. That's what we are trying to capture with
the notion of \emph{contraction in mean square}.

\begin{defn}
  The random system ${\bf x}_{i+1} = {\bf f}_i({\bf x}_i,i,\xii_i(\omega))$,
  is said to be \emph{contracting in mean square}
  if there exists an uniformly positive definite metric $M_i({\bf
    x},i,\xii_i)$ such that:
  $$
  \exists 0 \leq \eta < 1  \quad
  \expect( \| \delta {\bf z}(i,\xii_i) \|^2 ) \leq
  \|\delta {\bf z}_o\|^2 \eta^i
  $$
\end{defn}

\begin{thm}\label{thm:expectDiscret}

  Given the random system , $ {\bf x}_{i+1} = {\bf f}_i({\bf
    x}_i,i,\xii_i(\omega))$
  note $\sigma_{f_i}({\bf x},i,\xii_i)$ the largest singular value of the
  generalized Jacobian of ${\bf f}_i$ at time i 
  according to some metric $M$. \\
  
  A sufficient condition for the system to be \emph{contracting in
    mean square} is that 
  \begin{itemize}
  \item
    the random process $\sigma_{f_i}({\bf x},i,\xii_i)$
    can be bounded independently from ${\bf x}$, ie there exists a
    stochastic process $\eta_i$ such that  
    $$
    \forall {\bf x} \quad \sigma_{f_i}({\bf x},i,\xii_i) \leq
    \eta_i \qquad \as
    $$
  \item
    the stochastic process $\eta_i$ is constituted of
    independent random variables 
  \item
    the expectation of the random variables $\eta_i^2$ can
    be uniformly bounded by some $0 \leq \eta < 1$
    $$
    \forall i \quad \expect(\eta_i^2) \leq \eta
    $$
  \end{itemize}
\end{thm}

\begin{pf*}{Proof.}
  Note $F_i$ the discrete generalized Jacobian of $f_i$,
  we have again 
  $$
  {\bf \delta {\bf z}}_{i+1}^T {\bf \delta {\bf z}}_{i+1}
  \le  \sigma^2_{f_i}({\bf x}_i,i,\xii_i)
  \ {\bf \delta {\bf z}}_i^T {\bf \delta {\bf z}}_i
  \leq \eta_i^2 \ {\bf \delta {\bf z}}_i^T {\bf \delta {\bf z}}_i
  \qquad \as
  $$
  Introducing the expectation value of $\| \delta {\bf z} \|^2$, we
  use that independence between $X$ and $Y$ is defined as the
  uncorrelation of $f(X)$ and $g(Y)$ for all mesurable functions $f$
  anf $g$.
  $$
  \expect(\| {\bf \delta {\bf z}}_{i+1} \|^2 )
  \ \leq \ \expect(\eta_i^2 \| {\bf \delta {\bf z}}_i \|^2) =
  \expect(\eta_i^2) \expect(\| {\bf \delta {\bf z}}_i \|^2) 
  $$
  and hence
  \begin{equation*}
    \expect( \| \delta {\bf z}_i \|^2) \le \| \delta {\bf z}_o
    \|^2 \eta^i
  \end{equation*}
\qed
\end{pf*}

\subsection{Stochastic gradient}

Let us have a look to a stochastic way of minimizing a function,
highly used in computational neuroscience community, called
\emph{stochastic gradient}. 

The idea is to use the traditional minimization by 
\emph{gradient descent} but we want to avoid explicit computation of
the gradient as it is generally of high cost or even infeasible. 
For that, we introduced a stochastic perturbation which has the
role of a ``random differentiation''.

Let $\hat{P} = P + \Pi$ be the perturbation of
the state $P$ with respect the vector $\Pi$ of 
stochastic processes $\Pi_i$ .
Define the discrete system
$$
P_{n+1} = P_n - \mu.(\mathcal{E}(\hat{P_n}) -\mathcal{E}(P_n))\Pi
$$  
with $\mu > 0$.

Providing that the $\Pi_i$'s are mutually uncorrelated and of
auto-correlation $\sigma^2$ (ie. 
$\mathbb{E}(\Pi_i.\Pi_j) = \sigma^2.\delta_ {i,j}$),
the system satisfies:
$$
\delta P_{n+1} = \delta P_{n} - \mu.
(\sum_k \frac{\partial^2 \mathcal{E}}{\partial p_k \partial p_i}(P^*) 
\Pi_k.\Pi_j)_{i,j}.\delta P_{n}
$$
where $P^*$ is given by the finite difference theorem.
So, by taking the expectation:
$$
\expect(\delta P_{n+1}) = (I - \mu.\sigma^2.\frac{\partial^2 \mathcal{E}}
{\partial P_n^2}).\expect(\delta P_{n})
$$

So the system is contracting in mean square if

\begin{itemize}
\item
  $\frac{\partial^2 \mathcal{E}}{\partial P_n^2} > 0$ that is $\mathcal{E}$
  is strictly convex.
\item
  $\mu.\sigma^2\frac{\partial^2 \mathcal{E}}{\partial P_n^2} < I$ 
  (that is $\mu.\sigma^2$ sufficiently small)
\end{itemize}

\subsection{Nonlinear random system : the continuous time case}
\label{sct:continuous} 

We have seen in subsection \ref{sct:discrete} that the notion of
contraction for discrete-time varying systems harmonizes well with
stochastic calculus. Unfortunately the story is less straightforward
in the continuous time case. Nevertheless, as outlined in
introduction, for some practical reasons, we can restrict our intention to the
case of random differential systems as define in \cite{arnold98}. Let
us briefly summarize the technical background.

We want to define the stochastic extension of deterministic
differential systems as fellows.   

\begin{equation} \label{eq:continuous} 
  \dot{\bf x}_t = {\bf f}({\bf x}_t,t,\xii_t(\omega))
\end{equation}

where $\xii_t$ is a continuous stochastic process and
${\bf f}$ is a sufficiently smooth function, namely continuously
differentiable with respect to ${\bf x}_t$ (the condition on ${\bf f}$
can be reduced to a lipschitz condition, but as we need
differentiability in the rest of this paper, we prefer to assume it right
now).

But this formulation does not make sense for every kind of continuous
processes. Typically, when dealing with White noise process, the
right-hand part of equation \ref{eq:continuous} does not present
finite variation. 
In order to overcome this difficulty, we will assume that
$\xii_t$ is a ``nice'' stationary stochastic process whose 
trajectories are cadlag (for the french ``continue à droite et avec
des limites à gauche''), ie. are right continuous with left-hand
limits.

Arnold proved in \cite{arnold98} that under some assumption on 
${\bf f}$, equation \ref{eq:continuous} admits a unique solution which
is a global flow, whereas in general it is just a local flow. 

\begin{thm}[\cite{arnold98}]
  Suppose that $\xii$ is cadlag and ${\bf f} \in \mathcal{C}^1$.
  Then, equation \ref{eq:continuous} admit a unique maximal solution
  which is a local random dynamical system continuously differentiable.
  If furthermore ${\bf f}$ satisfies 
  $$
  \| {\bf f}(\xii,{\bf x})\| \leq \alpha(\xii) \|{\bf x}\| + \beta(\xii)
  $$
  where $t \rightarrow \alpha(\xii_t(\omega))$ and $t \rightarrow
  \beta(\xii_t(\omega))$ are locally integrable, then the solution is
  a global RDS.
\end{thm}

Thus, we cannot assume that random differential system defines a
unique continuous trajectory for every $\omega$. 
This problem is also known in the deterministic case where Vidyasagar has
shown the prevalence of differential equations despite our knowledge
of only very restrictive characterization. Indeed, the set of
equations admitting a unique solution is non-meager, whereas the set
of equations we are able to exhibit is meager. 
That is, ``practically all'' equations admit a unique solution whereas
we can characterize ``practically none'' of them !
So we will assume in the rest of this paper that the solution of the
differential equation exists and is a unique continuously
differentiable RDS.

All those restriction are rather technical and we refer the
interesting reader to \cite{arnold98} for further explanations. We can
yet give two slogans reformulating intuitions coming from those
restrictions :
\begin{description}
\item[The perturbation is memoryless]
  The noise appearing in the right-hand side of equation
  \ref{eq:continuous} is \emph{memoryless} in the sense that only the
  value of the perturbation at time $t$ enters into the generator
  ${\bf f}$.
\item[The perturbation do have small variations]
  The cadlag condition is a nice way to avoid problems generated by
  dramatically varying processes like White noise process while allowing
  interesting discontinuous processes such as jump Markov processes.
\end{description}

From now on, when we will talk about \emph{random differential
  system}, we assume that the solution of equation
\ref{eq:continuous} exists and is a unique continuously differentiable
RDS. We also suppose that all the processes we are dealing with are
stationary and have cadlag trajectories. 

\subsubsection{Almost sure contraction}

We can now define the notion of almost sure contraction for the
continuous-time case.

\begin{defn}

  A random differential system $\dot{\bf x}_t = {\bf f}({\bf x}_t,t,{\bf
    \xii}_t(\omega))$ is said to be \emph{almost surely contracting}
  if there exists an uniformly positive definite metric $ M({\bf
    x},t,\xii_t)$ and $\eta < 0$ such that:
  $$
  \mathbb{P}\{\omega \in \Omega, \ \limsup_{t \rightarrow \infty}
  \frac{1}{t} \log(\| \delta {\bf z} \|) \leq \eta \} = 1
  $$ ie. the difference between two trajectories tends almost surely 
  to zero exponentially.
\end{defn}

We now state conditions for a system to be almost surely
contracting. The traditional contraction analysis requires that the
largest eigenvalue of the general Jacobian is uniformly bounded by a negative
constant. The stochastic version of it mainly requires that the
largest eigenvalue which define a process, is bounded by a
process which follows the law of large number and of expectation
uniformly negative.

\begin{thm}\label{thm:asc}
  Given the system equations , $ \dot{\bf x}_t = {\bf f}({\bf x}_t,t,{\bf
    \xii}_t(\omega))$ note $\lambda_f({\bf x},t,\xii_t)$
  the largest eigenvalue of the generalized Jacobian of ${\bf f}$ at time t 
  according to some metric $M$.
  
  A sufficient condition for the system to be \emph{almost surely
  contracting} is that 
  \begin{itemize}
  \item
    the random process $\lambda_f({\bf x},t,\xii_t)$
    can be bounded independently from ${\bf x}$, ie 
    there exists a stochastic process $\eta_t$ such that 
    $$\forall {\bf x} \quad \lambda_f({\bf x},t,\xii_t) \leq
    \eta_t \qquad \as
    $$
  \item
    the stochastic process $\eta_t$ follows the strong
    law of large number
    $$ 
    \frac{1}{t} \int_0^t \eta_t \rightarrow
    \frac{1}{t} \int_0^t \expect(\eta_t)
    \qquad \as
    $$
  \item
    We can uniformly bound the expectation of the $\eta_t$
    with some $\eta < 0$
    $$\forall t \quad \expect(\eta_t) \leq \eta
    $$
  \end{itemize}
\end{thm}

\begin{pf*}{Proof.}
  We make a strong use of the basic idea of the original proof.
  $$
  \frac{1}{2} \ \frac{d}{dt} ({\bf \delta {\bf z}}^T {\bf \delta {\bf z}}) 
  \le 
  \lambda_f({\bf x}_t,t,\xii_t)
  \ {\bf \delta {\bf z}}^T {\bf \delta {\bf z}}
  \qquad \as
  $$
  and hence
  $$
    \| \delta {\bf z} \| \le \| \delta {\bf z}_o \| \ e^{\
      \int\limits_o^t \eta_t dt}
    \qquad \as
  $$
  
  Since $\eta_t$ verifies the law of large numbers, we have almost
  surely 
  $$
  \begin{array}{rcl}
    \limsup \frac{1}{t} \log \| \delta {\bf z}
    \| & \leq & \limsup \frac{1}{t} \int\limits_o^t \eta_t \ dt \\
    & = &\limsup \frac{1}{t} \int\limits_o^t
    \expect(\eta_t) dt \\
    & \leq & \eta
  \end{array}
  $$
  \qed
\end{pf*}

\begin{rem}
  It is reassuring that if we take a continuous random system satisfying
  conditions above, then the ``discrete envelop'' defined by $X_{n+1} =
  \exp [\int_{P_n} \lambda_f({\bf x}_t,t,\xii_t)] X_{n}$ is a
  discrete almost surely contracting system.
\end{rem}

\subsubsection{Contraction in mean square}

As it is the case in discrete-time case, we would like to find
sufficient that guarantee the contraction in mean square of our system
\begin{equation} \label{eq:continuous_average}
\expect(\| \delta {\bf z}(t,{\xii_t}) \|^2) \leq 
\|\delta {\bf z}_o\|^2 e^{\eta  t}
\end{equation}
Unfortunately, we have seen that this property required discrete independent
stochastic processes, whose continuous counterpart are processes like the
White noise process. As we have refused to deal with that kind of processes,
we need to find a stronger condition that yet ensure a similar
constraint on the average trajectory.

That's why we are moving to coarse-grained version of equation
\ref{eq:continuous_average}, namely where the property is guaranteed only
for a discrete sample of the average trajectory. This property will be
assessed when dealing with stochastic process which are coarse-grain
independent, as define below.

\begin{defn}
  A random differential system $\dot{\bf x}_t = {\bf f}({\bf
    x}_t,t,\xii_t(\omega))$ is \emph{coarse-grain contracting} if
  there exists a metric $ M({\bf x},t,\xii_t)$ and a
  partition $t_1 < t_2 < \ldots$ such that: 
  $$
  \exists 0 \leq \eta < 1 \ , \ \forall i \quad \expect(\| \delta {\bf
    z}_{t_i} \|^2) \leq \|\delta {\bf z}_o\|^2 \eta^i
  $$
\end{defn}

To guarantee this property, we have to deal with particular kind of
continuous stochastic processes which satisfy a condition of
independence in a coarse-grain scale.

\begin{defn}[coarse-grain independence]
  A continuous stochastic process $\eta_t$ is said to be \emph{coarse
    grain independent with respect to a partition} $(P_i)$ of $\Preal$ if
  $$
  \mbox{the random variables } 
  \eta_i = \int_{P_i} \eta_t \ dt
  \mbox{ are independent}
  $$
\end{defn}

\begin{rem}
\begin{enumerate}
\item
  By a partition $(P_n)$, we mean equivalently a strictly increasing
  infinite sequence $t_1 < t_2 < \ldots$ or a sequence of intervals $P_0 =
  [0,t_1], P_i = (t_i,t_{i+1}]$ for $i\geq1$. 
\item
  In case of Gaussian or uniform random variables, the condition is
  satisfied if two random variables lying in two different sets of the
  partition are always independent.
\end{enumerate}
\end{rem}

\paragraph*{Example of coarse grain independent process}
We will now define the typical type of coarse grain independent
process we have in mind. Take a partition $(P_n)$ of
$\Preal$ and an independent stochastic process $G_n(\omega)$. 
  
Define the process $\gamma_t(\omega) = G_n(\omega)$ for $t \in P_n$.

Then $\gamma_t(\omega)$ is a coarse grain independent process. In that
case, each trajectory is piecewise constant and we have 
$$
 \int_{P_n} \gamma_{t} dt = |P_n|.G_n \mbox{ define an independent
   stochastic process}
$$

\begin{figure}[h] \label{fig:coarse-grain}
  \begin{center}

    \includegraphics[width=\textwidth]{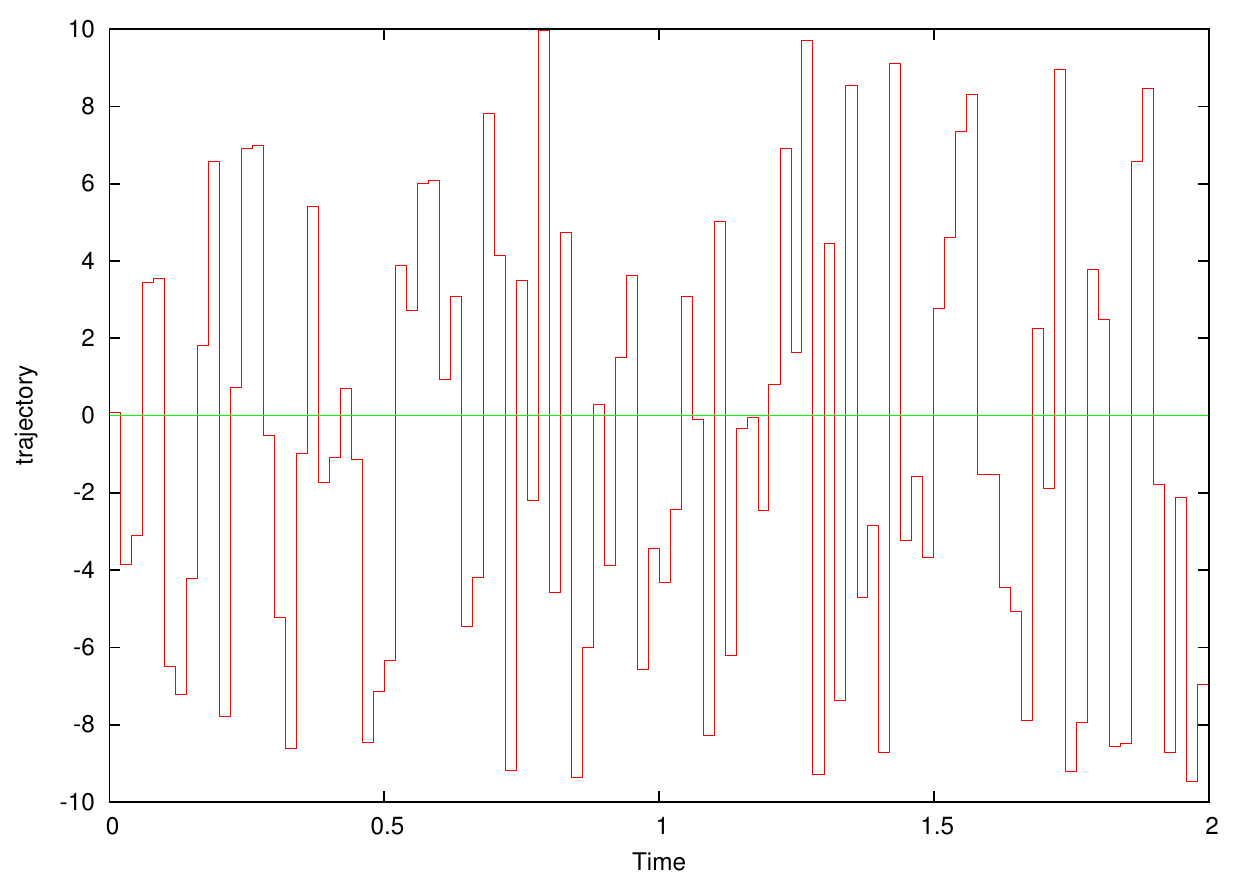}
    \caption{Example of a trajectory of a coarse grain independent
      process} 

  \end{center}
\end{figure}

\begin{thm}\label{thm:expect}
  Given the system equations , $ \dot{\bf x}_t = {\bf f}({\bf x}_t,t,{\bf
    \xii}_t(\omega))$ note $\lambda_f({\bf x},t,\xii_t)$
  the largest eigenvalue of the generalized Jacobian of ${\bf f}$ at time t 
  according to some metric $M$.
  
  A sufficient condition for the system to be \emph{coarse-grain
    contracting} is that
  \begin{itemize}
  \item
    the random process $\lambda_f({\bf x},t,\xii_t)$
    can be bounded independently from ${\bf x}$, ie 
    there exists a stochastic process $\eta_t$ such that 
    $$
    \forall {\bf x} \quad \lambda_f({\bf x}_t,t,\xii_t) \leq
    \eta_t
    \qquad \as
    $$
  \item
    the process $\eta_t$ is a coarse-grain
    independent stochastic process with respect to a partition
    $(P_n)$. 
  \item
    We can uniformly bound the expectation of the $e^{\int_{P_n}
      \eta_t}$ with some $0\leq \eta < 1$ 
    $$\forall n \quad
    \expect((e^{\int_{P_n} \eta_t})^2) \leq \eta
    $$
  \end{itemize}
\end{thm}

\begin{pf*}{Proof.}
   $$
     \frac{1}{2} \ \frac{d}{dt} ({\bf \delta {\bf z}}^T {\bf \delta {\bf z}}) 
     \le 
     \lambda_f({\bf x}_t,t,\xii_t)
     \ {\bf \delta {\bf z}}^T {\bf \delta {\bf z}}
     \le \eta_t
     \ \delta {\bf z}^T \delta {\bf z}
     \qquad \as
  $$
  Which leads to
  $$
  \|\delta {\bf z}\| \leq
  \| \delta {\bf z}_0 \| \ e^{\int_0^t \eta_t dt}
  \qquad \as
  $$
  Thus, we can define the system $Z_{n+1} = e^{\int_{t_n}^{t_{n+1}}
    \eta_t dt} Z_n$ and $Z_0 = \| \delta {\bf z}_0 \|$, which
  satisfies 
  $$
  \|\delta {\bf z}_{t_n}\| \leq  Z_n
  $$
  By definition of coarse-grain stochastic process and as
  $\expect((e^{\int_{P_n} \eta_t})^2) \leq \eta$, we can applied
  theorem \ref{thm:expectDiscret} to conclude on the contraction in
  mean square of $Z_n$ with rate $\eta$
  $$
  \expect (\| {\bf \delta {\bf z}}_{t_i} \|^2)
  \leq
  \expect(Z_i^2) \leq \| \delta {\bf z}_0 \|^2 \eta^i
  $$
\qed
\end{pf*}

\begin{rem}
  The condition imposed on the process $\eta_t$, namely
  $\expect((e^{\int_{P_n} \eta_t})^2) \leq \eta < 1$, is really
  different from the condition we have seen for the almost sure
  contraction $\expect(\eta_t) \leq \eta < 0$.
\end{rem}

\section{Part 2 : Noise's dependency}

Let us now turn on a very special case of perturbed systems, namely
when the impact of the noise does not depend on the state space.

\begin{prop}
  Consider a contracting system $ \dot{\bf \widehat{x}} = {\bf f}({\bf
    \widehat{x}},t)$ and
  take a perturbed version of it $ \dot{\bf x}_t = {\bf f}({\bf x}_t,t) +
  {\bf   \xii}_t(\omega)$.
  The system is automatically both contracting on average and almost
  surely contracting.
\end{prop}

The proof is obvious as it is the case mention above of a system which
is contracting for every $\omega$. Let us now study the mean and the
variance of the unique solution.

\subsection{Study of the average trajectory} \label{sct:average trajectory}

Consider a contracting system $ \dot{\bf \widehat{x}} = {\bf f}({\bf
  \widehat{x}},t)$ in the metric $M = \mathnormal{I}$,
take a perturbed version of it $ \dot{\bf x}_t = {\bf f}({\bf x}_t,t) +
{\bf   \xii}_t(\omega)$. 
Then, assuming that 
\begin{itemize}
\item
  $\expect(\xii_t) = 0$
\item
  $\forall t \ \|\xii_t\| < \alpha$ almost surely with
  $\expect(\alpha) < \infty$ 
\end{itemize}
we have that 
$$
\forall t \quad \expect({\bf x}_t) = {\bf \widehat{x}}_t  
$$

\begin{pf*}{Proof.}
  Let ${\bf x_\Delta} = {\bf \widehat{x}} - \expect({\bf x}_t)$ and
  consider
  \begin{eqnarray*}
    {\bf x}_\Delta^T \frac{d}{dt}({\bf x}_\Delta) & = &
    {\bf x}_\Delta^T ({\bf f}({\bf \widehat{x}},t) - \expect({\bf
    f}({\bf x}_t,t))) \quad \mbox{by dominated convergence theorem}\\
    & = & {\bf x}_\Delta^T \expect({\bf f}({\bf \widehat{x}},t) -
    {\bf f}({\bf x}_t,t)) \\ 
    & = & {\bf x}_\Delta^T \expect(\int_0^1 \frac{\partial f}{\partial
   {\bf x}} ({\bf \widehat{x}} + c {\bf x}_\Delta,t) dc \ {\bf
   x_\Delta}) \\
    & \leq & \expect(\int_0^1 \lambda_f^{max}{\bf x}_\Delta^T {\bf
    x_\Delta}) \\
    & = & \lambda_f^{max} \| {\bf x}_\Delta \| ^2
  \end{eqnarray*}
  So we have that $\| {\bf x}_\Delta \| \leq \| {\bf
    x}_{\Delta_0} \| e^{2 \lambda_f^{max} t}$.
  But $\| {\bf x}_{\Delta_0} \| = 0$.
\qed
\end{pf*}

\subsection{Study of deviation}

\begin{thm}
  Take random system $ \dot{\bf x}_t = {\bf f}({\bf x}_t,t) +
  {\xii_t(\omega)}$ satisfying the conditions of the subsection
  \ref{sct:average trajectory}. Suppose now that the deviations of the
  $\xii_t$ are uniformly bounded
  $$
  \expect(\|\xii_t\|) \leq \sigma
  $$
  The deviation of ${\bf x}_t$ is then majored by the maximum
  deviation $\sigma$ in the following way:
  $$
  \expect(\| {\bf x_1} - {\bf x_2} \|) - \expect(\| {\bf {\bf x_1} - {\bf x_2}}
  \|_0) e^{\lambda_f^{max} \ t} \leq
  \frac{2 \sigma}{|\lambda_f^{max}|} (1 - e^{\lambda_f^{max} \ t}))
  $$ 
\end{thm}

\begin{pf*}{Proof.}
 Let ${\bf \widetilde{x}} = {\bf x_1} - {\bf x_2}$
 Let us look at $\|{\bf \widetilde{x}}\|$
 $$
 \left.
 \begin{array}{l}
   \dot{\bf x_1} = {\bf f}({\bf x_1},t) +
       {\xii_1} \\
       \dot{\bf x_2} = {\bf f}({\bf x_2},t) +
       {\xii_2} \\
 \end{array}
 \right\}
 \Rightarrow
 \frac{d}{dt}{\bf \widetilde{x}} = \int_0^1 \frac{\partial f}{\partial
 {\bf x}} ({\bf x_2} + c {\bf \widetilde{x}},t) dc \ {\bf
 \widetilde{x}} + {\xii_1} - \xii_2
 $$
 Multiply by ${\bf \widetilde{x}^T}$, it becomes:
 $$
 {\bf \widetilde{x}^T} \frac{d}{dt}{\bf \widetilde{x}} = 
 {\bf \widetilde{x}^T} \int_0^1 \frac{\partial f}{\partial
 {\bf x}} ({\bf x_2} + c {\bf \widetilde{x}},t) dc \ {\bf
   \widetilde{x}} + {\bf \widetilde{x}^T} (\xii_1 - \xii_2)
 $$
 $$
 \frac{1}{2} \frac{d}{dt}{\| \bf \widetilde{x} \| ^2} \leq 
 \lambda_f^{max} \|{\bf \widetilde{x}}\|^2 + \|{\bf \widetilde{x}}\| \
 \|\xii_1 - \xii_2\|
 $$
 So, by dominated convergence theorem again
 $$
 \frac{d}{dt} \expect(\|{ \bf \widetilde{x}} \|) \leq
 \lambda_f^{max} \expect(\|{\bf \widetilde{x}}\|) +
 \expect(\|\xii_1 - \xii_2\|) \leq
 \lambda_f^{max} \expect(\|{\bf \widetilde{x}}\|) +
 2 \sigma
 $$
 Solving $\frac{d}{dt} \expect(\|{ \bf \widetilde{x}} \|) =
 \lambda_f^{max} \expect(\|{\bf \widetilde{x}}\|) +
 2 \sigma$ and using
 the positiveness of all terms in the equation (which means that
 replacing $=$ by $\leq$ just make the slope of increasement smaller),
 we have :
 $$
 \expect(\| {\bf \widetilde{x}} \|) - \expect(\| {\bf \widetilde{x}}
 \|_0) e^{\lambda_f^{max} \ t} \leq
 \frac{2 \sigma}{|\lambda_f^{max}|} (1 - e^{\lambda_f^{max} \ t}))
 $$
 \qed
\end{pf*}

\subsection{Oscillator Synchronization}

Consider two identical Van der Pol oscillators couples as
$$
\ddot{\bf x}_1 + \alpha({\bf x}_1^2 - 1)\dot{\bf x}_1 + w^2 {\bf
  x}_1 = \alpha {\epsilon_1}(\dot{\bf x}_2 - \dot{\bf x}_1)
$$
$$
\ddot{\bf x}_2 + \alpha({\bf x}_2^2 - 1)\dot{\bf x}_2 + w^2 {\bf
  x}_2 = \alpha {\epsilon_2} (\dot{\bf x}_1 - \dot{\bf x}_2)
$$

where 
\begin{itemize}
\item 
  $\alpha > 0$, $w >0$
\item
  {$\epsilon_1$}, {$\epsilon_2$} are stationary processes
\end{itemize}

Using (Combescot, Slotine 2000),
we can show that $x_1 \xrightarrow[t \rightarrow \infty]{} x_2$ when
$$
\mathbb{E}(\epsilon_1) + \mathbb{E}(\epsilon_2) > 1
$$

\begin{figure}[h]
  \begin{center}

    \includegraphics[width=0.49\textwidth]{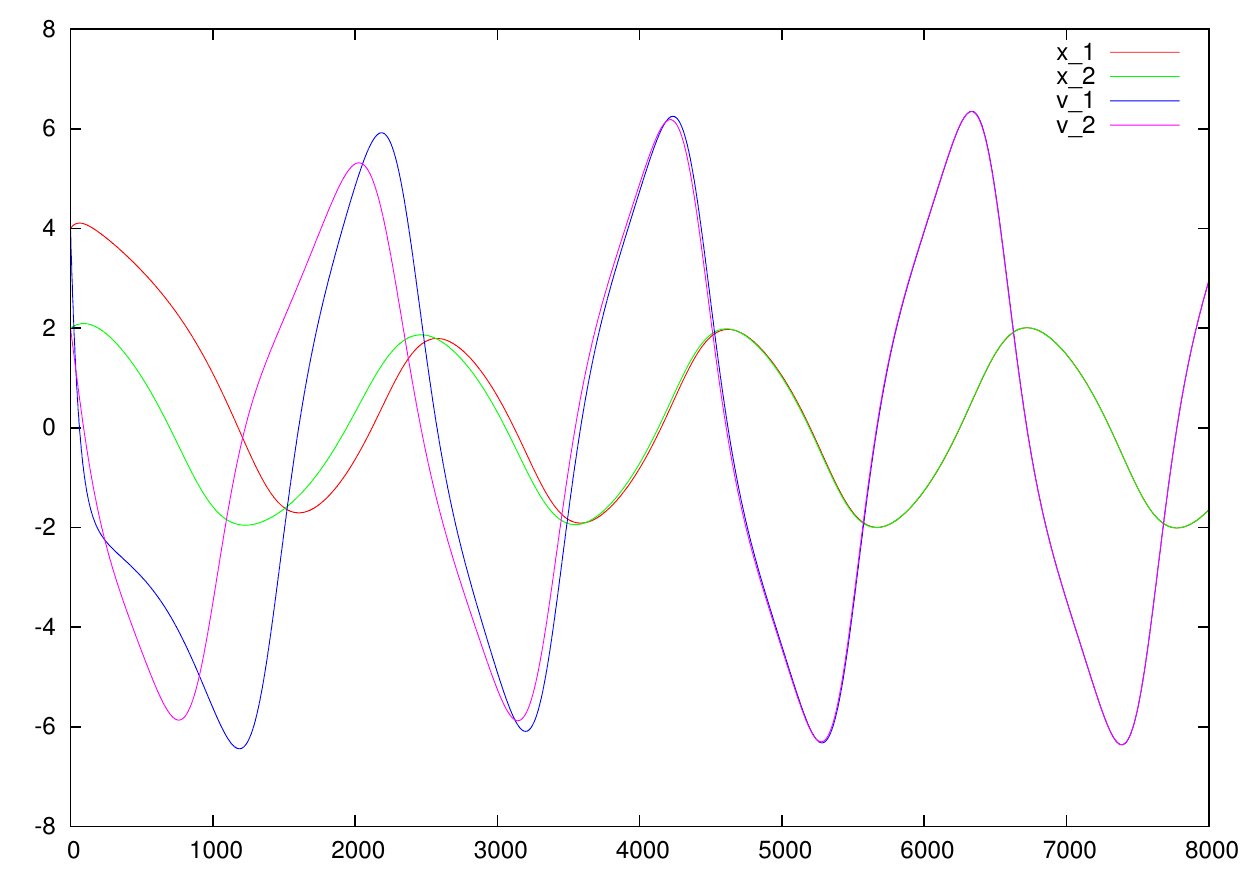}
    \includegraphics[width=0.49\textwidth]{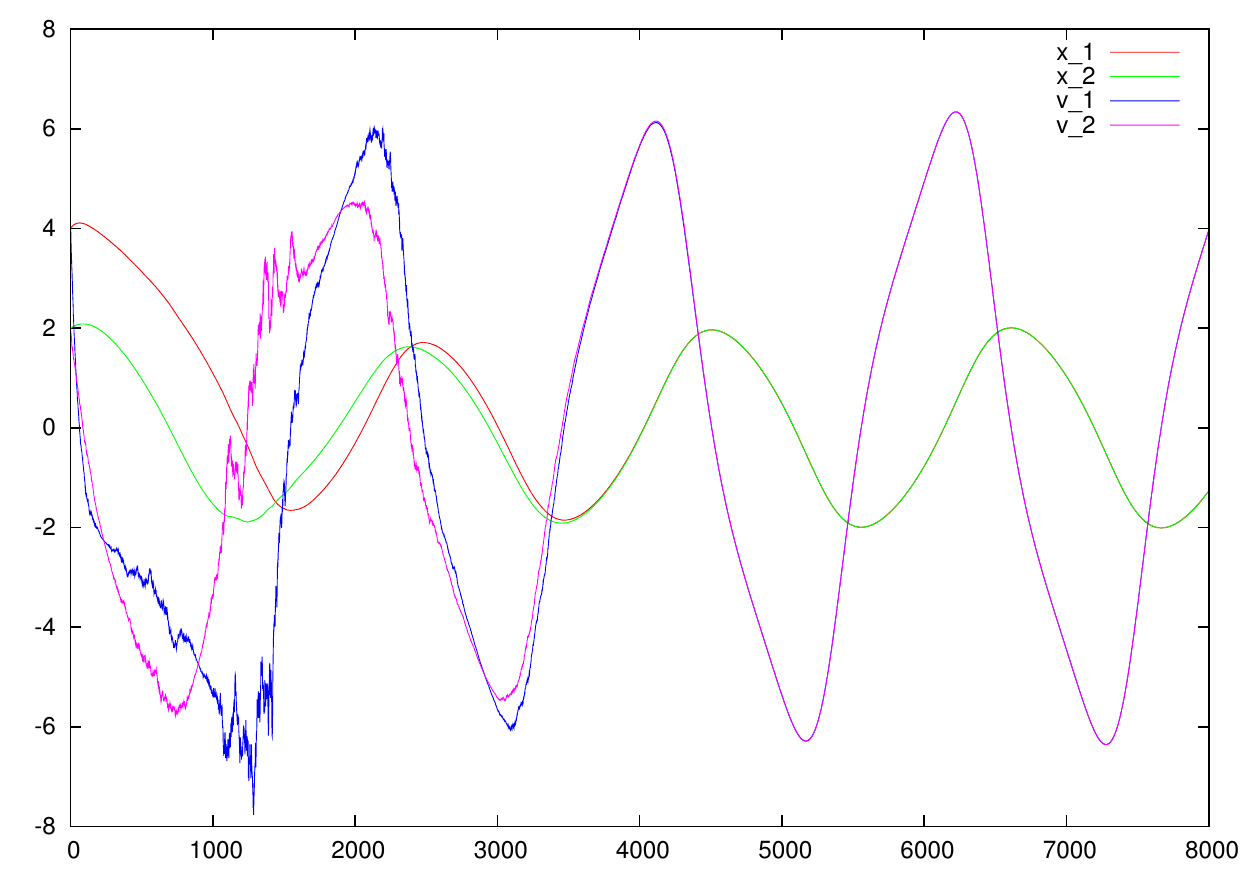}

    \caption{A comparison between stochastic synchronization of two
      Van der Pol oscillator with and without noise ($\epsilon_1,\epsilon_2 \in [-40,40]$)} 

  \end{center}
\end{figure}

\begin{figure}[h] \label{fig:inputNoise}
  \begin{center}
    
    \includegraphics[width=0.49\textwidth]{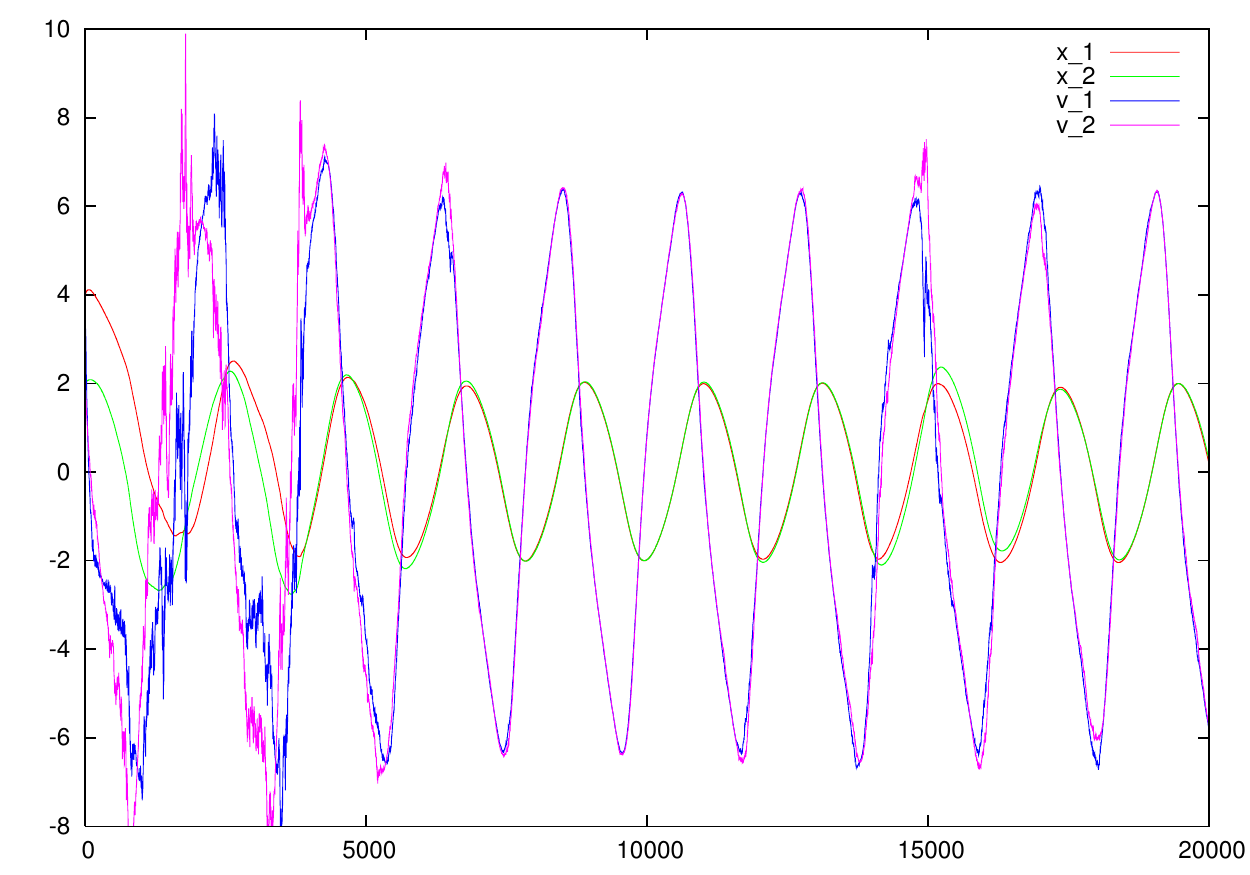}
    \includegraphics[width=0.49\textwidth]{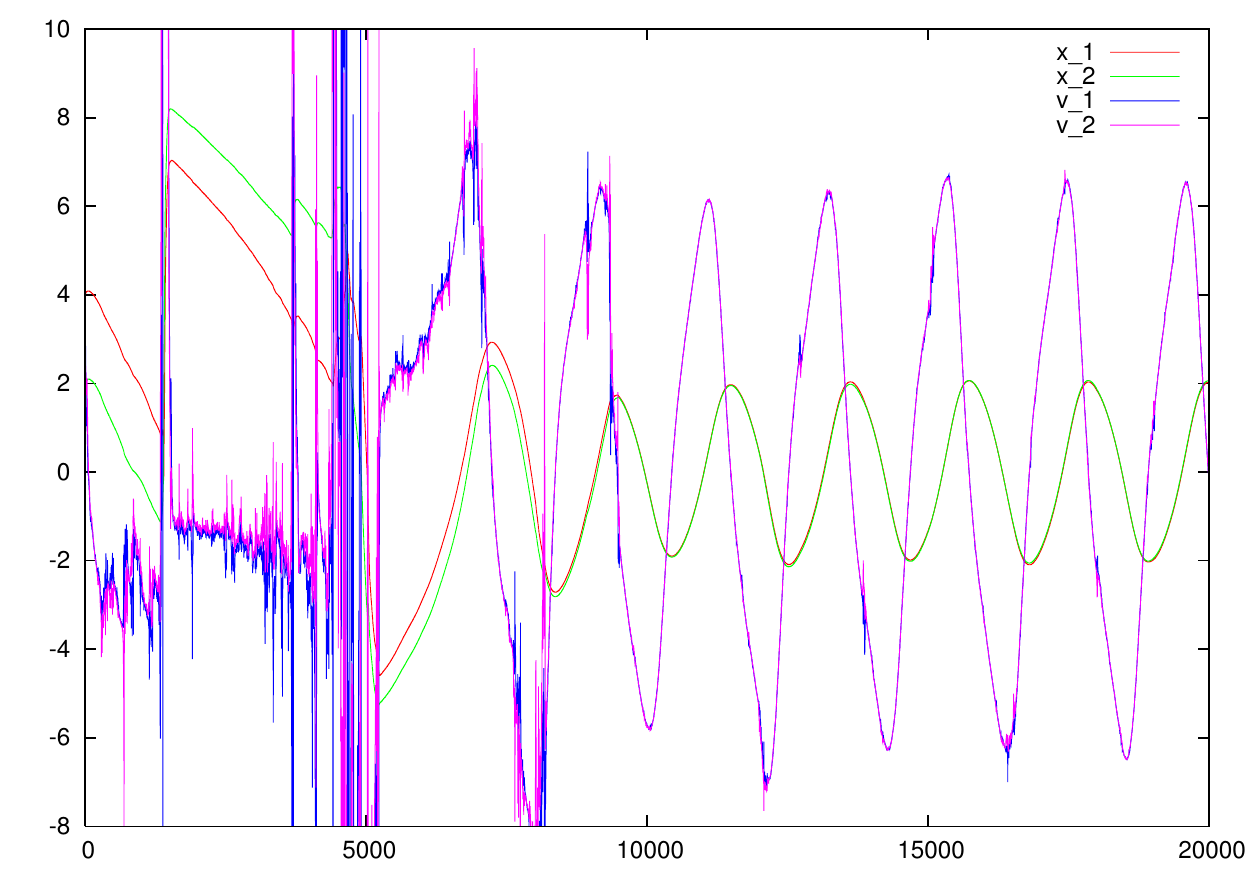}
    \caption{Stochastic synchronization with noise in input. Left part
    obtained with $\epsilon_1,\epsilon_2 \in [-80,80]$, right part with
    $\epsilon_1,\epsilon_2 \in [-400,400]$}  
    
  \end{center}
\end{figure}

Remark that we can add noise in the input of both oscillators, the
synchronization still occurs on average (fig. \ref{fig:inputNoise}).  

\bibliographystyle{apalike}
\bibliography{biblio.bib}


\end{document}